\documentclass{autart}               

\usepackage{amsmath,amsfonts,framed,enumitem}
\usepackage{algorithm,algpseudocode}
\newcommand{\subscript}[2]{$#1 _ #2$}

\usepackage{graphicx}          
\usepackage[dvips]{epsfig}    
\usepackage{comment,xcolor}      
\usepackage{mathrsfs,placeins}
                              
\allowdisplaybreaks

\begin{document}

\begin{frontmatter}

\title{Stable and Robust LQR Design via Scenario Approach} 

\thanks[footnoteinfo]{This paper was not presented at any IFAC 
meeting. Corresponding author: Anna Scampicchio.}

\author[Padova]{Anna Scampicchio}\ead{anna.scampicchio@phd.unipd.it},    
\author[UW]{Aleksandr Aravkin}\ead{saravkin@uw.edu},
\author[Padova]{Gianluigi Pillonetto}\ead{giapi@dei.unipd.it}   

\address[Padova]{University of Padova, Padova, Italy.}  
\address[UW]{University of Washington, Seattle, WA}        


\begin{abstract}              Linear Quadratic Regulator (LQR) design is one of the most classical optimal control problems, whose well-known solution is an input sequence expressed as a state-feedback.
In this work, finite-horizon and discrete-time LQR is solved under stability constraints and uncertain system dynamics. The resulting feedback controller balances cost value and closed-loop stability. 
Robustness of the solution is modeled using the scenario approach, without requiring any probabilistic description of the uncertainty in the system matrices. 
The new methods are tested and compared on the Leslie growth model, where we control population size while minimizing a suitable finite-horizon cost function.     
\end{abstract}

\end{frontmatter}

\section{Introduction}

Optimal control is a mature branch of applied mathematics, arising from the interaction of calculus of variations \cite{pontryagin}, optimization \cite{LuenbergerVectorSpaceOpt}, differential geometry \cite{BressanPiccoli} and control system design \cite{lewis}. This area has recently attracted the machine learning community thanks to the parallels with reinforcement learning \cite{recht}: both disciplines deal with the problem of finding an optimal input that minimizes some cost function under the constraint of the (possibly uncertain) system dynamics. Moreover, the same problem is the main ingredient in Model Predictive Control (MPC) \cite{rawlings2009model}. In the stochastic setting, optimal control sees its dual in the theory of optimal state estimation \cite{todorovconf}. 

The best known problem in optimal control is likely the Linear Quadratic Regulator (LQR) design, where the cost is taken as quadratic and the differential/difference equation governing the system is assumed to be linear and deterministic \cite{Anderson}. The assumptions above yield a \textit{convex} problem. It turns out that the optimal input is a state-feedback, whose expression is related to the solution of a Riccati equation \cite{lancaster1995algebraic}, which will be Difference/Differential (DRE) for finite-horizon, and Algebraic (ARE) for infinite-horizon setups. Usually, generic optimal control problems do not allow for a closed form solution.

This work focuses on discrete-time, finite-horizon LQR problem in the following context:
\begin{enumerate}
\item casting the solution to be a state-feedback, the final closed-loop system has to be stable;
\item system matrices are affected by uncertainty, available through samples rather than a functional probabilistic description. 
\end{enumerate}

In the optimal control literature, stability has been thoroughly investigated only in the infinite-horizon setting \cite{Kalman60contributionsto}. Indeed, as pointed out in \cite{Bitmead1991RiccatiDA}, discussing stability in the classic finite-horizon LQR solution is not reasonable, since the DRE has no finite escape properties: at each updating step, all terms in the DRE are well defined and finite. 
In \cite{martenssonrantzer}, a finite-horizon LQR problem is studied but stability results are achieved by interpreting it as an approximation of an infinite-horizon problem. 
Optimizing a finite-horizon cost is the fundamental step in MPC, where the issue of stability is carefully studied, e.g see  \cite{denicolaomagnistabilitymmpc}, \cite{Mayne2000ConstrainedMP}. 
MPC is an iterative method that solves a finite-horizon optimal control problem at each time step; next, only the first input is applied to the system, and then the procedure is repeated on the updated state. 
The algorithm resembles dynamic programming approaches \cite{Bertsekas}, but instead of looking for a closed-loop solution, MPC seeks an open-loop input sequence that can be computed online. For this reason, it avoids the so-called ``curse of dimensionality" and has less computational burden arising from constraints on inputs and states. Since MPC solution is expressed as open-loop, stability questions have to be posed focusing on the orbits induced by the designed inputs. Originally, the discussion focused only on the horizon length, but in more recent work  two main approaches concerning stability have emerged (\cite{rawlings2009model}, \cite{grune2017nonlinear} \cite{Mayne2000ConstrainedMP}):
\begin{enumerate}[label=(\alph*)]
\item treat the cost as a Lyapunov function, carefully tuning its parameters;\label{enum:stab1}
\item force the final state to lie in the neighbourhood of the equilibrium point (usually the origin). \label{enum:stab2}
\end{enumerate}
Approach \ref{enum:stab1} was suggested by \cite{Kalman60contributionsto} for the infinite-horizon LQR, while \ref{enum:stab2} was first used by \cite{kleinman}, \cite{thomas} and \cite{kwon}.
In \cite{BEMPORAD19981255}, both approaches are used, along with an extra set of variables that lets the original system be asymptotically stable.

\paragraph*{Contribution}
The main novelty of this work is to fill the gap between dynamic programming and open-loop solutions, searching for a static feedback matrix that yields a stable closed-loop system by keeping a finite-horizon cost as low as possible. Hence, the problem is formulated in terms of trade-off between stability and cost value. Moreover, the approaches we propose to deal with stability do not rely on \ref{enum:stab1} or \ref{enum:stab2}: specifically, we do not tune the parameters in the cost in order to obtain a suitable Lyapunov function, and we do not impose a dead-beat behaviour (as e.g. in \cite{Calafiore2013RobustMP}) or any constraint on the value of the final state. Our point of view leads to nonconvex problems that are anyway effective in terms of stability, cost value and computational time. For all of them, we will perform the convergence analysis for the optimization problems involved.\\

{The second issue}, i.e. \textit{uncertainty in system matrices}, prevents the application of classic methods such as dynamic programming dealing with uncertain parameters \cite{Bertsekas}, \cite{hou2015dynamic}: in fact, no probabilistic description of system dynamics uncertainty is available. We complement our new optimal control technique with the so-called \textit{scenario approach}, a data-driven method introduced in \cite{calafiorecampi} that exploits just samples of the uncertain quantities. The scenario approach was also used for model predictive control applications \cite{mpc}, \cite{Calafiore2013RobustMP}, but not in the form explored in this paper.\\ 
  
The new robust approach to LQR is especially useful when the original system is unstable with uncertain parameters and the aim is to stabilize it in such a way that a cost function over a finite time horizon is minimized. Such a situation arises in many cases e.g. to control the growh of a cell population. One of the possible models describing this dynamics in discrete time is the one introduced by Leslie in \cite{leslie}. It will be used in the numerical experiments.\\

\paragraph*{Roadmap}
This work proceeds as follows. In Section \ref{sec:prelim}, we review concepts in linear dynamical systems and scenario approach optimization. Section \ref{sec:classic} introduces the LQR problem in the discrete-time, finite-horizon case, providing its formulation in four equivalent ways that are all used in developing the new approaches. In particular, one viewpoint is generalized to solve the unconstrained case with uncertainty, while the others are used as starting points to develop the approaches dealing with stability constraints. The latter are investigated in Section \ref{sec:stability}, both in the deterministic and in the robust case. Section \ref{sec:numexp} gathers all the numerical tests. Finally, in Section \ref{sec:conclusions}, conclusions are drawn.

\paragraph*{Notation} Throughout this paper, $\mathcal{L}(\cdot)$ will be used to denote all the Lagrangians that appear, and similarly $\lambda_t$ is denoting a general multiplier at time $t$. A gradient of a function $g$ w.r.t. a vector $v$ will be written as $\bigtriangledown_v g$. As regards matrix manipulations, $\otimes$ is used for Kronecker products, $vec(\cdot)$ for the vectorization operator and $\mathscr{K}_{d_1d_2}$ indicates the $d_1d_2 \times d_1d_2$ commutation matrix. Moreover, $\text{diag}(\cdots)$ denotes a (block-)diagonal matrix, $1_N$ is a column vector composed by ones, $I_N$ is the $N \times N$ identity matrix and $0_{n,m}$ is a $n \times m$ matrix of zeroes. If a matrix $P$ is (semi)definite positive, then it will be expressed as $P\succ 0 \;(P \succeq 0)$. Eigenvalues are denoted by the letter $\sigma$. Finally, the Frobenius norm of a matrix will be defined as $\| \cdot \|_F$, while $\| \cdot \|$ denotes the standard 2-norm.

\section{Preliminaries}\label{sec:prelim}
\subsection{Review on Linear Systems Theory}
We start with the deterministic setup. 
We consider dynamical systems that are linear, time-invariant, and with dynamics described by the difference equation
\begin{equation}\label{eq:dyn}
x_{t+1} = Fx_t + Gu_t, \qquad
x_0 = \bar{x},
\end{equation}
where, for all $t$, $x_t \in \mathbb{R}^n$ and $u_t \in \mathbb{R}^m$ are the state and the  input respectively. We  use $x$ and $u$ to denote the column "meta-vectors" obtainined by stacking the state/input sequences, e.g. $u = [u_0^{\top} \; u_1^{\top} \; \cdots]^{\top}$. \\
Considering the solution of the homogeneous (i.e. inputless) case, the dynamics are governed by the eigenvalues of $F$. In particular, if their modulus is strictly less than or bigger than 1, a mode is asymptotically stable or unstable (we do not discuss the case in which $|\cdot | = 1$). A famous criterion to assess if a system is asymptotically stable without computing its eigenvalues is due to Lyapunov \cite{1961iii}. 
Going back to the original system \eqref{eq:dyn}, the main task in control theory consists in designing the input $u$ in order to avoid unstable state trajectories. This is effectively done by imposing the input to be a state-feedback ruled by some gain matrix $K$, i.e. $u_t = Kx_t$: in this way, the dynamics becomes $x_{t+1} = (F + GK)x_t$. Matrix $K$ can be suitably designed to stabilize the system if the system is controllable: hence we will assume throughout the paper that each system has this property. We refer to \cite{Anderson} for a review on the notions of reachability/controllability and observability/detectability. Under the hypothesis of controllability, we are ready to state an extended version of Lyapunov's Theorem above mentioned, with both the classical statement for linear systems and its version as a Linear Matrix Inequality (LMI) \cite{DEOLIVEIRA1999},\cite{finsler1937uber} \cite{BEFB:94}.
\begin{thm}\label{lyp}
The following are equivalent:
\begin{itemize}
\item system \eqref{eq:dyn} is asymptotically stable with $u_t=Kx_t$
\item $\forall\, \Xi \succ 0,\; \exists P \succ 0$ such that $(F+GK)^{\top}P(F+GK) - P = -\Xi$
\item for some matrices $P$, $C$ and $D$ of suitable dimensions,
\begin{equation}\label{eq:stabLMI}
\begin{bmatrix}
P & FC + GD\\
C^{\top}F^{\top} + D^{\top}G^{\top} & C + C^{\top} - P
\end{bmatrix} \succ 0.
\end{equation}
Moreover, knowing matrices $C$ and $D$, it results that $K = DC^{-1}$ whenever $C$ is invertible.
\end{itemize}
\end{thm}
To obtain a robust solution, the region of uncertainty is parametrized by $\delta \in \Delta$. While in theory $\Delta$ may be infinite, 
we assume it is endowed with a $\sigma-$algebra $\mathcal{D}$ and a probability \textsf{Prob}. Uncertainty will enter system matrix $F = F(\delta)$. In this work we assume that matrix $G$ is deterministic; however, the case in which $G = G(\delta)$ does not require any methodological novelty, and the tools below introduced can be easily extended to deal with this situation. 

\subsection{On the scenario approach}
Consider the following \textit{convex} optimization problem
\begin{equation}\label{eq:originalMinMax}
 \min_{\theta \in \Theta} \max_{\delta \in \Delta} c(\theta,\delta),
\end{equation}
which can be rephrased in epigraphic form \cite{Boyd}:
\begin{equation}\label{eq:scenarioIntermediate}
\min_{\alpha, \theta} \quad \alpha \quad \text{s.t. } c(\theta,\delta) \leq \alpha \quad \forall \, \delta \in \Delta.
\end{equation}
Since $\Delta$ is uncountable, this problem involves an infinite number of constraints. The so-called scenario approach is a data-driven relaxation of this problem. Instead of considering a standard worst-case approach, which may be too conservative, optimization is carried out w.r.t. $N$ of samples $\lbrace \delta^{(i)} \rbrace_{i=1}^N$, thus considering
\begin{equation}\label{eq:scenario}
\min_{\alpha, \theta} \quad \alpha \quad \text{s.t. } c(\theta,\delta^{(i)}) \leq \alpha \quad \forall \, \delta^{(i)}, \, i=1...N.
\end{equation}
The power of the theory of the scenario approach in the convex setup lies in the fact that it is possible to know in advance the number $N$ of samples in order to guarantee a certain level of robustness with desired confidence. In fact, define the following quantities:
\begin{itemize}
\item $V(\theta) = \textsf{Prob}\{\delta \in \Delta: c(\theta,\delta)-\alpha > 0 \}$ is the probability that a solution $\hat{\theta}$ of problem \eqref{eq:scenarioIntermediate} violates the constraint $c(\theta,\delta) - \alpha \leq 0$.
\item The parameter $\epsilon \in (0,1)$ is aimed at defining the concept of $\epsilon-$level solution, which is: $\theta \in \Theta$ such that $V(\theta) \leq \epsilon$.
\item $\beta \in (0,1)$ is the risk of failure, i.e. $\beta$ bounds the risk that our solution is not an $\epsilon-$level solution.
\end{itemize}
Then, the following Theorem holds \cite{CampiGarattiM03}:
\begin{thm}\label{thm:scenario}
Fix $\epsilon$ and $\beta$ above defined, and let $d$ be the number of optimization variables. If the number of scenarios $N$ satisfies \begin{equation}\label{eq:Nsc}
N \geq \frac{2}{\epsilon}\Big(\ln\frac{1}{\beta} + d - 1 \Big)
\end{equation} 
then, with probability no smaller than $1 - \beta$, either problem \eqref{eq:scenario} is infeasible (then also the original robust convex program is); or, if it is feasible, then the optimal solution $\hat{\theta}$ is $\epsilon-$robustly feasible.
\end{thm}

This strong result depends on a consequence of Helly's Theorem \cite{Helly1923}: in particular, the number of \textit{support constraints} (i.e. constraints that, if removed, improve the solution) never exceeds the number of optimization variables if the problem is convex. 

\paragraph*{What if the problem is non-convex?}
The result given above cannot be applied if the cost is nonconvex. If this is the case, then no guarantee about the number of support constraints can be provided: hence, a solution of the minmax problem can be given without any a-priori guarantee on its optimality/robustness.  Nevertheless, significative results are present in \cite{alamotempocamacho}, where the Vapnik-Chervonenkis theory \cite{Vapnik1998} is exploited; then in \cite{grammaticononconv}, where the nonconvex support is reduced to a convex setup in order to apply Theorem \ref{thm:scenario}; finally, in \cite{campinonconv}, where the less conservative and most flexible approach is presented. In particular, the following Theorem holds:
\begin{thm}\label{thm:nonconvScenario}
Consider problem \eqref{eq:scenario}. Fix confidence parameter $\beta \in (0,1)$ and a certain number $N$ of scenarios. Consider a function $\epsilon: \{0, \cdots, N \} \rightarrow [0,1]$ such that
\begin{displaymath}
\epsilon(N)=1, \qquad \sum_{k=0}^N {N\choose k}(1 - \epsilon(k))^{N-k}=\beta.
\end{displaymath}
Then it holds that
\begin{displaymath}
\textsf{Prob}\Big(V(\bar{\theta}_N > \epsilon(\bar{s}_N))\Big) \leq \beta
\end{displaymath}
where $\bar{s}_N$ is the cardinality of the considered support subsample and $\bar{\theta}_N$ is the corresponding optimal solution.
\end{thm}
Therefore, when dealing with nonconvex problems, the procedure consists in selecting the support constraints, solving the problem and then retrieving \textit{a posteriori} the robustness of the solution according to Theorem \ref{thm:nonconvScenario}.

\section{LQR and its solution in the scenario approach}\label{sec:classic}
\subsection{General theory for LQR}
The classic discrete-time, finite-horizon LQR problem is formulated as follows:
\begin{align}\label{pb}
\min_{u}{x_T^{\top}Sx_T + \sum_{t=0}^{T-1}(x_t^{\top}Qx_t + u_t^{\top}Ru_t)} \\
s.t. \begin{cases} x_{t+1}=Fx_t + Gu_t \nonumber\\  x_0 = \bar{x} 
\end{cases}
\end{align}
where $Q, \,S \succeq 0$ and $R \succ 0$.
The problem can be posed/solved in four equivalent ways, denoted by (P1), (P2), (P3) and (P4) below described. 
All of these viewpoints are used in new developments later  in the paper.

\subsubsection{Approach (P1): use the dynamics to rewrite the problem only in terms of the inputs}
Use the constraint on the system dynamics \eqref{eq:dyn} to rewrite the objective \eqref{sec:classic} in terms of $u$ only, obtaining the problem 
\begin{equation}
\label{P1:QP}
\min_u u^{\top}Bu + 2a^{\top}u
\end{equation}
whose solution is $u = -B^{-1}a$. 
The expressions for $a$ and $B$ are obtained by substituting in the objective the expression derived from \eqref{eq:dyn} for each $t$:
\begin{displaymath}
x_t = F^{t}\bar{x} + \sum_{s=0}^{t-1}F^{t-1-s}Gu_s.
\end{displaymath}
Then, the following properties for multiple sums are used:
\begin{displaymath}
\sum_{t=1}^{T-1}\sum_{k=1}^t = \sum_{k=1}^{T-1}\sum_{t=k}^{T-1}\text{ and }
\sum_{t=1}^{T-1}\sum_{i=1}^{t}\sum_{j=1}^{t} = \sum_{i=1}^{T-1}\sum_{j=1}^{T-1}\sum_{t=\max(i,j)}^{T-1}.
\end{displaymath}
The final expressions are two block matrices, whose $(i,j)$ block is expressed as 
\begin{align*}
[B]_{i,j=1...T} &= \Big[G^{\top}(F^{T-i})^{\top}SF^{T-j}G\Big]_{i,j=1...T}+[R]_{i=j}+\\
  & + \Bigg[\sum_{t=max(i,j)}^{T-1}G^{\top}(F^{T-i})^{\top}QF^{t-j}G\Bigg]_{i,j=1...T-1}\\  
  [a]_{i=1...T} &= \Big[\bar{x}^{\top}{(F^T)}^{\top}SF^{T-i}G\Big]_{i=1...T} +\\ &+\Bigg[\sum_{t=i}^{T-1}\bar{x}^{\top}(F^t)^{\top}QF^{t-i}G\Bigg]_{i=1...T-1}
\end{align*}
\subsubsection{Approach (P2): matricial form and KKT conditions}
Consider objective and constraints written in matricial form as
\begin{equation}
\label{P2:QP}
\min_u x^{\top}\bar{Q}x + u^{\top}\bar{R}u \quad \text{s.t. } A_1x + A_2u - b = 0
\end{equation}
where
$\bar{Q} = diag(I_T\otimes Q, S)$, $\bar{R}=I_T \otimes R$,
\begin{equation}\label{eq:relationsP2} 
A_1 = \begin{bmatrix}
I \\
-F & I \\
0 &-F &I\\
0 &\cdots &-F &I
\end{bmatrix} { \scriptstyle \in \mathbb{R}^{n(T+1)\times n(T+1)}},\end{equation}
$$A_2 = \begin{bmatrix}
0 &\cdots &0\\
-G\\
0 &\cdots &-G
\end{bmatrix} { \scriptstyle \in \mathbb{R}^{n(T+1)\times mT}},\, b = \begin{bmatrix}
\bar{x}\\ 0 \\ \vdots
\end{bmatrix} {\scriptstyle \in \mathbb{R}^{n(T+1)}}. $$
Then, compute the KKT conditions \cite{Boyd}, \cite{kuhn1951} on the Lagrangian $\mathcal{L}(x,u,\lambda)=x^{\top}\bar{Q}x + u^{\top}\bar{R}u + \lambda^{\top}(A_1x + A_2u - b)$: 
\begin{equation}\label{U2optcond}
\begin{cases}
\bigtriangledown_x \mathcal{L} = 0 \rightarrow 2\bar{Q}x + A_1^{\top}\lambda = 0\\
\bigtriangledown_u \mathcal{L} = 0 \rightarrow 2\bar{R}u + A_2^{\top}\lambda = 0 \\
\bigtriangledown_{\lambda} \mathcal{L} = 0 \rightarrow A_1x + A_2u - b=0
\end{cases}
\end{equation}
Since the problem is convex, necessary optimality conditions are also sufficient for the solution. We obtain thus \begin{align*} u &= (\bar{R} + A_2^{\top}A_1^{-\top}\bar{Q}A_1^{-1}A_2)^{-1}A_2^{\top}A_1^{-\top}\bar{Q}A_1^{-1}b\\ &= \bar{R}^{-1}A_2^{\top}A_1^{-\top}\bar{Q}x. \end{align*}
\subsubsection{Approach (P3): Objective manipulation and DRE} The objective can be rewritten as follows \cite{kwakernaak}:
\begin{align*}
&x_0^{\top}M_0x_0 + x^{\top}_T[S - M_T]x_T +\\ 
\sum_{t=0}^{T-1} \Big([R + &G^{\top}M_{t+1}G]u_t + G^{\top}M_{t+1}Fx_t\Big)^{\top}\times \\ 
&\times [R + G^{\top}M_{t+1}G]^{-1}\times \\ 
\times\Big([R &+ G^{\top}M_{t+1}G]u_t + G^{\top}M_{t+1}Fx_t \Big).
\end{align*}
This form is valid for any arbitrary matrix sequence $\lbrace M_t \rbrace_{t=0}^T$, but it can be proved via e.g. Bellman's optimality principle \cite{Bellman} that the optimal sequence has to satisfy the DRE
\begin{equation}\label{dre}
\begin{aligned}
 M_t & = Q + F^{\top}M_{t+1}F \\
 & - F^{\top}M_{t+1} G(R + G^{\top}M_{t+1}G)^{-1}G^{\top}M_{t+1}F.
 \end{aligned}
\end{equation}
Then, by inspection of the objective, it results that the DRE has to satisfy the boundary condition $M_T = S$ and that $u_t = \underbrace{-(R + G^{\top}M_{t+1}G)^{-1}G^{\top}M_{t+1}F}_{K_t}x_t$ for each $t$.

\subsubsection{Approach (P4): Pontryagin Maximum Principle}
Let us recall the main result \cite{pontryagin}: 
\begin{thm}[Pontryagin's maximum principle]
Consider the problem
\begin{displaymath}
\min_{u,x} {\psi(x_T) + \sum_{t=0}^{T-1}l(x_t,u_t)} \; \text{s.t. } x_{t+1}= f(x_t,u_t),\;x_0=\bar{x}.
\end{displaymath}
and consider its optimal solution. Define the Hamiltonian $$H(x,u,\lambda)=l(x,u) + \lambda^{\top}(f(x,u)).$$ Then there exists a sequence $\lbrace \lambda_t \rbrace_{t=0}^{T}$ such that  $u_t$, $x_t$ and $\lambda_t$ satisfy the following conditions:
\begin{displaymath}
\begin{cases}
x_{t+1} = f(x_t, u_t)\\
x_0 = \bar{x}\\
\lambda_t = \frac{\partial l}{\partial x_t} l(x_t, u_t) - f(x_t,u_t)^{\top}\lambda_{t+1}\\
\lambda_T = -\frac{\partial \psi}{\partial x_t}(x_T)\\
\frac{\partial H(x_t, u_t, \lambda_{t+1})}{\partial u} = 0
\end{cases}
\end{displaymath}
\end{thm}
The proof of this Theorem is made exploiting variational considerations  \cite{luenberger1979introduction} and does not invoke KKT conditions.
In the LQR setup, $\psi(x_T)=x_T^{\top}Sx_T$, $l(x_t,u_t)= x_t^{\top}Qx_t + u_t^{\top}Ru_t$ and the optimality conditions listed in Pontryagin's maximum principle reduce to this set of equations:
\begin{equation}\label{pontryagin}
\begin{cases}
x_{t+1}=Fx_t + Gu_t\\
x_0 = \bar{x}\\
\lambda_t = -2Qx_t + F^{\top}\lambda_{t+1}\\
\lambda_T = -2Sx_T\\
2Ru_t - G^{\top}\lambda_{t+1} = 0.
\end{cases}
\end{equation}
From this, we obtain the control law $u_t = \frac{1}{2}R^{-1}G^{\top}\lambda_{t+1}$. By construction, $x_t$ and $\lambda_t$ depend linearly on $x_0$, thus hinting that there is a linear relation between $\lambda_t$ and $x_t$: we  then guess that $\lambda_t = -2P_tx_t$ for some matrix sequence $P_t$. This operation is referred to as "sweep method" in \cite{brysonho}. We then obtain $u_t = -(R + G^{\top}P_{t+1}G)^{-1}G^{\top}P_{t+1}Fx_t$. To find out which rule has to be satisfied by $P_t$, we consider the dynamics of the multiplier $\lambda_t$:
\begin{align*}
\lambda_t &= -2Qx_t + F^{\top}\lambda_{t+1}\\ \rightarrow {2}P_t x_t &= {2}Qx_t + {2}F^{\top}P_{t+1}(Fx_t + Gu_t)
\end{align*}
where, substituting the expression for $u_t$ and dropping the dependence from $x_t$, we recognize the Riccati equation \eqref{dre}.
\subsubsection*{Comments}
We have shown four equivalent approaches that can be used to obtain the solution for the classic LQR problem. 
The optimal input can be expressed both in closed form (P1, P2) and as a state-feedback (P2, P3, P4).\\
From a computational point of view, the solution obtained via (P3) and (P4) is computed in $O(T)$ time, since it basically requires only the computation of the DRE. 
In contrast, naive implementation of (P1) and (P2) requires inversion of full-state matrices and thus needs $O(T^3)$ time. 
The equivalence allows us to think of the efficient methods $(P3)-(P4)$ as numerical subroutines to solve the 
full-state optimization problems~\eqref{P1:QP} and~\eqref{P2:QP}, which is very convenient when developing extensions to the approach.

\subsection{Extension to the scenario case}
The robust unconstrained LQR problem is solved drawing inspiration from approach (P1), whose formulation aligns with \eqref{eq:scenario}. Each sample of system matrix $F(\delta^{(i)})$ yields different matrices $B^{(i)}$ and $a^{(i)}$, so that the problem becomes the convex program
\begin{equation}
\min_{u,\alpha} \alpha \quad \text{s.t. } u^{\top}B^{(i)}u + 2{a^{(i)}}^{\top}u \leq \alpha \;\forall i=1...N,
\end{equation}
where $N$ can be chosen according to Theorem \ref{thm:scenario}.
The solution can be implemented e.g. in \textsc{Matlab} using \texttt{cvx}. A possible way to assess the correctness of the solution is to compare the dual variables that can be given as output in \texttt{cvx} with the optimality conditions. In particular, define $\bar{B} = \text{diag}(B^{(1)},\cdots B^{(N)})$ and $\bar{a}= [a^{(1)} \cdots a^{(N)}]$ that yield $J(u) = (I_N \otimes u^{\top})\bar{B}(1_N \otimes u) + 2\bar{a}^{\top}u$. Notice that $J(u)$ is just a compact way to express the $N$ constraints. Then, taking the Lagrangian $\mathcal{L}(\alpha,u,\mu) = \alpha + \lambda^{\top}(\alpha1_N - J(u))$ and computing the gradients w.r.t. $\alpha$ and $\mu$, we get
\begin{displaymath}
1_N^{\top}\lambda +1 = 0, \quad \alpha 1_N - J(u) \geq 0.
\end{displaymath}

\section{Adding the stability constraint}\label{sec:stability}
We first consider the deterministic setup. In Section \ref{sec:classic} we saw that the optimal solution is described by a sequence of feedback matrices $\{K_t\}_{t=0}^{T-1}$, which in turn is derived by the DRE ``initialized" at $M_T = S$. In our setup, we imagine applying  this sequence to the original system in order to minimize the cost, and then let the system evolve in closed loop with the last feedback $K_{T-1} = -(R + G^{\top}SG)^{-1}G^{\top}SF$. \textit{Our aim is to obtain a final closed-loop system that is asymptotically stable}. Depending on the choice of $S$ and $R$,  $F + GK_{T-1}$ may be unstable. So, instead of considering $S$ and/or $R$ as additional optimization variables, we look for a solution which is expressed in terms of a \textit{constant feedback} able to stabilize the system. This is going to be a suboptimal solution w.r.t. the one obtained with (P1) --- (P4): the advantage of the LQ setup is that it will be possible to compare how far is the stabilizing solution from the classic unconstrained one. \\
The new problem in the deterministic setup has the following structure:
\begin{align}\label{eq:stabprob}
\min_{u,x} x_T^{\top}Sx_t + \sum_{t=0}^{T-1}u_t^{\top}Ru_t + x_t^{\top}Qx_t\\
\text{subject to} \nonumber
\begin{cases}
u_t = Kx_t\\
(F + GK) \text{ is stable}\\
x_{t+1} = Fx_t+ Gu_t\\
x_0 = \bar{x}
\end{cases}
\end{align}
We consider two possible strategies to handle the stability constraint:
\begin{enumerate}
\item include it in the set of constraints as LMI, drawing inspiration from Theorem \ref{lyp};
\item exploit the theory for the infinite-horizon case, aiming at finding an approximate solution of the ARE. The Riccati equation, under suitable conditions on the matrices involved, is guaranteed to yield a stabilizing solution.
\end{enumerate}

Now we present those approaches in more detail. The first  produces a method that we denote by (S0). The second can be treated in different ways, and  yields (S1), (S2) and (S$_{\infty}$).


\subsection{First strategy: Lyapunov stability included as constraint - approach (S0)}\label{sec:S0}
Taking the problem expressed as \eqref{eq:stabprob}, substituting the expression for $u_t$ in the objective and including the stability constraint along the lines of Theorem \ref{lyp}, we obtain the following problem:
\begin{align*}
&\min_{x, K,P,C,D} x_T^{\top}Sx_T + \sum_{t=0}^{T-1}x_t^{\top}(Q +  K^{\top}RK)x_t \\
&\text{s.t. } \begin{cases}
x_{t+1}=(F+GK)x_t, \qquad x_0 = \bar{x}\\
\begin{bmatrix}
P & FC + GD \\ (FC + GD)^{\top} & C + C^{\top} - P
\end{bmatrix} \succeq \xi I_{2n}\\
KC=D.
\end{cases} \nonumber
\end{align*}
The LMI above differs from the one in \eqref{eq:stabLMI}: indeed, semidefinite programming allows only for non-strict inequalities, and our aim is to impose the constraint matrix to be positive definite. We therefore choose the parameter $\xi$ to be small and positive (e.g. $10^{-5}$).\\
The first crucial step is to relax the problem by inserting the constraint $KC=D$ (see Theorem \ref{lyp}) as a regularization term in the objective. 
We then  include the linear constraints via a Lagrange formulation: 
\begin{equation}
\label{eq:S0obj}
\begin{aligned}
 \bar{\mathcal{L}}(x,\lambda,K) := &x_T^{\top}Sx_T + \lambda_0^{\top}(x_0 - \bar{x}) +\\ + &\sum_{t=0}^{T-1} x_t^{\top}(Q + K^{\top}RK)x_t  \\
 + &\lambda_{t+1}(x_{t+1} - (F+GK)x_t) \\
=& x^{\top}\bar{B}(K)x + \lambda^{\top}(\bar{A}(K)x - b),
\end{aligned}
\end{equation}
where $\bar{B}(K) = \bar{Q} + \text{diag}\Big(I_{T}\otimes (K^{\top}RK), 0_{n,n}\Big)$ and $\bar{A}(K)$ has the structure of matrix $A_1$ in \eqref{eq:relationsP2}, with $F$ replaced by $F + GK$. Moreover, introducing the shorthand notation 
\begin{equation}
E=
\begin{bmatrix}
\xi I_n - P & -(FC+GD)\\ -(FC+GD)^{\top} & \xi I_n + P - C - C^{\top}
\end{bmatrix}, \nonumber
\end{equation}
the problem can be written as 
\begin{equation}
\min_{x,K,\atop P,C,D}\max_{\lambda,\Gamma} \bar{\mathcal{L}}(x,\lambda,K) + \frac{1}{2\mu}\| KC-D\|_F^2 + tr(\Gamma^{\top}E).
\end{equation}
The overall problem is nonconvex; nevertheless, it is a bi-convex program in the variables ($x$, $K$) and $(P,C,D)$. We solve it using 
alternating minimization~\cite{nonconvOptMachineLearning}, \cite{TsengAltMin}; the overall procedure is summarized in Algorithm \ref{alg:1}.

\begin{algorithm}
\caption{\normalsize \textbf{- Solution of (S0).} The inputs are system matrices $F$ and $G$, initial state $\bar{x}$, objective matrices $Q$, $R$ and $S$, parameters $\mu$ and $\xi$. The desired output is the (stabilizing) feedback matrix $K$.\label{algo:alt}}
\label{alg:1}
$i=0$;\\
initialize $K(0),\,C(0),\,D(0),\,P(0)$ randomly;\\
$i=1$;\\
\textbf{while} {not converge} \textbf{do}
\begin{enumerate}[label=\subscript{C}{{\arabic*}})]
\item predispose $v(\cdot)=\min_x\max_{\lambda}\bar{\mathcal{L}}(x,\lambda,\cdot)$ via \eqref{eq:C1optcond};
\item $\displaystyle K(i)= \min_K v(K) +\frac{1}{2\mu}\|KC(i-1)-D(i-1)\|^2_F$;
\item $\displaystyle P(i),\,C(i),\,D(i) = \min_{P,C,D}\|K(i)C-D\|^2_F$\\ subject to $\begin{bmatrix}
P & FC+ GD \\ (FC + GD)^{\top} & C + C^{\top} - P
\end{bmatrix} \succeq \xi I_{2n}$;
\end{enumerate}
\begin{itemize}
\item $i=i+1$;
\end{itemize}
\textbf{end while}
\end{algorithm}
\setlist[enumerate,1]{start=1}

The preliminary step, which is denoted by $(C_0)$, consists in partially minimizing over the primal and dual variables $x$ and $\lambda$, leaving a value function of $K$ alone:
\begin{equation}
\label{eq:value}
\begin{aligned}
v(K) &= \min_{x} \max_{\lambda} \bar{\mathcal{L}}(x, \lambda, K) \\
& = \bar{\mathcal{L}}(x(K), \lambda(K), K)  
\end{aligned}
\end{equation}
where the tuple $(x(K), \lambda(K))$ denotes the optimal primal-dual pair at a fixed $K$, and $\bar{\mathcal{L}}$ is as in~\eqref{eq:S0obj}.
Using problem structure, we have explicit expressions of these quantities:
\begin{equation}\label{eq:C1optcond}
\begin{cases}
x=\bar{A}(K)^{-1}b\\
\lambda = -\frac{1}{2}\bar{A}^{\top}(K)\bar{B}(K)\bar{A}(K)b.
\end{cases}
\end{equation}
At this point one should minimize the objective w.r.t. $K$: this step is denoted as $(C_1)$. The expressions $x(K)$ and $\lambda(K)$ depend on $K$ in a complex way, 
but this dependence does 
not affect the first derivative of the value function, in fact we do not need use the specific closed form expressions~\eqref{eq:C1optcond}
to implement the method; we need any computational routine that returns $(x(K), \lambda(K))$~\cite{aravkin2012estimating}. 
In other words, $\nabla v(K)$ can be completely captured using the values of $(x(K), \lambda(K))$, again through the structure~\eqref{eq:S0obj}:  
\begin{equation}
\label{eq:valuegrad}
\begin{aligned}
\nabla v(K) & = x(K)^{\top}\partial_K\bar{B}(K)x(K) \\
&+ \lambda(K)^{\top}(\partial_K\bar{A}(K))x(K) - b),
\end{aligned}
\end{equation}
where the differentials with respect to the matrix $K$ are written formally and must be correctly computed in coordinates. Since the derivative is now available, 
we can use L-BFGS to solve (C1). 
With some computations reported in the Appendix we obtain an explicit representation for the value function gradient:
\begin{equation}\label{eq:gradient}
\sum_{t=0}^{T-1}2RKx_tx_t^{\top} - G^{\top}\lambda_{t+1}x_t^{\top} + (KCC^{\top} - DC)/\mu.
\end{equation} 
The last subproblem is denoted as $(C_2)$ and consists in minimizing over $(P,C,D)$ while keeping $K$ fixed. We obtain the convex program 
\begin{align}\label{eq:C2}
&\min_{P,C,D} \| KC-D \|^2_F\\
\text{subject to }&\begin{bmatrix}
P & FC+ GD \\ (FC + GD)^{\top} & C + C^{\top} - P
\end{bmatrix} \succeq \xi I_{2n}. \nonumber
\end{align}
which is solved using \verb{cvx{ applying interior point methods \cite{NoceWrig06} to $\min_{P,C,D} \|KC-D\|_F^2 +  tr(\Gamma^{\top}E)$.
\setlist[enumerate,1]{start=0}

The block-alternating scheme summarized in Algorithm \ref{alg:1} consists in partially minimizing w.r.t. $x$, and then using 
an alternative update on the value function with respect to $K$ and $(P,C,D)$. 
The special structure of the problem ensures convergence according to the main results of ~\cite{tseng2001convergence}. In particular, we notice that the problem is structured as
\begin{equation}
\min_{z_1,z_2} g(z_1,z_2) \nonumber
\end{equation}
with $z_1=K$ and $z_2=(P,C,D)$. The value function $g$ is convex with respect to each variable block; Then, Theorem 4.1.a in \cite{tseng2001convergence} can be applied, and convergence to a stationary point is guaranteed.

\subsection{Second strategy: ARE theory - approaches (S$_{\infty}$), (S1) and (S2)}
This strategy draws inspiration for the results of the infinite-horizon LQR problem, whose objective is \eqref{pb} with $T \rightarrow +\infty$: \begin{displaymath}
\min_u \sum_{t=0}^{+\infty}x_t^{\top}Qx_t + u_t^{\top}Ru_t.
\end{displaymath}
Accordingly, its optimal solution is $u_t = -(R + G^{\top}M_{\infty}G)^{-1}G^{\top}M_{\infty}Fx_t$. Matrix $M_{\infty}$ comes from the asymptotic version of the DRE initialized at $0_{n,n}$, and is one of the possible solutions of the ARE \cite{lewis}. In particular, this important result holds:
\begin{thm}\label{thm:stabARE}
Let $(F,G)$ be controllable. $M_{\infty} \succeq 0$ is the unique stabilizing solution of ARE if and only if $(F,Q^{1/2})$ is detectable.
\end{thm}

Therefore, a suboptimal solution of the finite-horizon LQR problem can be expressed in terms of $M_{\infty}$, that is also stabilizing under the detectability condition. This approach will be denoted as (S$_\infty$). This solution is reliable only if we know that the detectability hypothesis in Theorem \ref{thm:stabARE} is verified for each scenario: if this is not the case, then numerical solvers fail in providing a solution. The other two approaches, denoted by (S1) and (S2), try to approximate the ARE solution. The three approaches are furtherly investigated below.

\subsubsection{Approach (S$_{\infty}$)}
The solution of the ARE is obtained numerically using methods that hinge on \cite{ArnoldLaub}. Under the hypotheses of Theorem \ref{thm:stabARE}, the matrix obtained coincides with the asymptotic solution of the DRE initialized at $0_{n,n}$. Moreover, for a sufficiently high $T$, the DRE solution aligns with $M_{\infty}$ after a transitory given by matrix $S$: in particular, if $S=M_{\infty}$, the DRE yields its stationary solution. \\

\subsubsection{Approach (S1)}
This approach draws inspiration from (P4). A key passage in the application of Pontryagin's Maximum Principle was setting $\lambda_t = -2P_tx_t$, and subsequent rearrangements of the optimality conditions yielded the solution to the LQR problem. Since a constant feedback matrix is searched, we approximate the optimal solution by setting $\lambda_t = 2L^{\top}Lx_t$ for some $L \in \mathbb{R}^{n\times n}$ and consider only the optimality condition $\frac{\partial H(x,u,\lambda)}{\partial u}$. Then we obtain $u_t = -(R + G^{\top}L^{\top}LG)^{-1}G^{\top}L^{\top}LFx_t$, that we then substitute in the original problem \eqref{pb}. We thus obtain the nonconvex problem
\begin{align}\label{eq:arefeedback}
&\min_L \bar{x}^{\top}(F^{\top} + K^{\top}G^{\top})^{T}S(F + GK)^{T}\bar{x} + \nonumber \\+ \sum_{t=0}^{T-1} &\bar{x}^{\top}(F^{\top} + K^{\top}G^{\top})^{t}(Q + K^{\top}RK)(F+GK)^t\bar{x} \nonumber\\
&\text{with } K = -(R + G^{\top}L^{\top}LG)^{-1}G^{\top}L^{\top}LF,
\end{align}
which is solved using a gradient-free method, e.g. Nelder-Mead's \cite{NoceWrig06} implemented with \texttt{fminsearch} in \textsc{Matlab}. 

\subsubsection{Approach (S2)}
This method is derived from (P3), where the objective had been rewritten using an (at first) arbitrary matrix sequence. Imposing that this sequence has to be constant, we obtain that the optimal feedback has the same structure of \eqref{eq:arefeedback}, and the objective is 
\begin{displaymath}
\min_L \bar{x}^{\top}L^{\top}L\bar{x} + \bar{x}^{\top}(F^{\top} + K^{\top}G^{\top})^{T}(S - L^{\top}L)(F + GK)^{T}\bar{x}.
\end{displaymath}
Again, this problem is solved via Nelder-Mead's method. We can see that $L^{\top}L=S$ is one of the solutions, and is the optimal one but is possibly non stabilizing. Nevertheless, the aim of this method (and of the previous one as well) is to try to find suboptimal solutions in which $\|x_t \|$ is decreasing as $t$ increases. 


\subsection{Extension with the scenario approach}
The previous part of Section \ref{sec:stability} was devoted to the deterministic setup only: now we want to consider the robust case. First of all, it can be noticed that all approaches exposed in Section \ref{sec:stability} yield \textit{nonconvex} problems: therefore, Theorem \ref{thm:scenario} cannot be applied. The only guarantee can be given a posteriori using Theorem \ref{thm:nonconvScenario}. The key passage consists in isolating an \textit{irreducible} support subsample, i.e. a set such that no element can be furtherly removed without changing the solution. A possible procedure, proposed in \cite{campinonconv}, consists in (1) solving the original problem with all $ N$ constraints, i.e. $\mathcal{S}= \{\delta^{(1)},...,\delta^{(N)}\}$, (2) for each $i=1...N$, consider $\mathcal{S}^{\prime} \leftarrow \mathcal{S}\backslash \delta^{(i)}$ and solve again with $\mathcal{S^{\prime}}$: if the result is unchanged, then $\delta^{(i)}$ is a support constraint. Having isolated a support subsample of cardinality $s_N^*$ with the procedure above sketched, we can proceed in evaluating $\epsilon(s_N^*)$ and get the measure of robustness of our solution. Following \cite{campinonconv}, function $\epsilon(\cdot)$ in Theorem \ref{thm:nonconvScenario} is chosen for all cases as 
\begin{equation}\label{eq:epsilon}
\epsilon(k)=\begin{cases} 
1 &\qquad \text{if }k=N\\
1 - \sqrt[N-k]{\frac{\beta}{N{N\choose k}}} &\qquad \text{otherwise.}\\
\end{cases}
\end{equation}
We now adapt all approaches to the robust case. To align with \eqref{eq:scenario} and to facilitate the support constraint selection procedure, all minmax formulations will be rephrased in epigraphic form. As proved e.g. in \cite{CaE:14}, this yields an equivalent statement of the problem.  

\textit{As regards (S0)}, Algorithm \ref{alg:1} is modified as follows. First, subproblem (C1) is expressed as
\begin{displaymath}
\min_K \max_{\delta^{(i)}} \bar{\mathcal{L}}\Big( x(\delta^{(i)},K), \lambda(\delta^{(i)},K), K \Big) + \frac{1}{2\mu}\|KC-D \|_F^2.
\end{displaymath}
At this point, subproblem (C2) presented in \eqref{eq:C2} becomes again a minmax problem, where the maximisation step is taken over all scenarios, i.e.
\begin{align}
&\qquad \quad \min_{P,C,D}\max_{\delta^{(i)}} \| KC-D \|^2_F \nonumber \\
\text{s. t. }&\begin{bmatrix}
P & F(\delta^{(i)})C+ GD \\ (F(\delta^{(i)})C + GD)^{\top} & C + C^{\top} - P
\end{bmatrix} \succeq \xi I_{2n}. \nonumber
\end{align}


\textit{As regards (S1) and (S2)}, let us call $J(\delta^{(i)}, L)$ the objective of both approaches. Then we have to solve for
\begin{displaymath}
\min_L \max_{\delta^{(i)}} J(\delta^{(i)}, L).
\end{displaymath} 
Having obtained the solution $L$, the feedback matrix $K = -(R + G^{\top}L^{\top}LG)^{-1}G^{\top}L^{\top}LF(\delta^{(i)})$ is built for every new scenario yielding $F(\delta^{(i)})$.\\
\textit{Lastly, for (S$_{\infty}$)}, we have to solve $N$ AREs numerically in order to find the scenario that maximizes the objective. If $(F(\delta^{(i)}),Q^{1/2})$ fails to be detectable for some $i$, no numerical solution of the ARE is provided. We deal with this situation by sampling $N$ scenarios, and then discard the ones that are not detectable. Hence, when computing $\epsilon(\cdot)$ of Theorem \ref{thm:nonconvScenario}, the starting number of scenarios is the one corresponding to the detectable ones.

\section{Numerical experiments}\label{sec:numexp}
In this Section, we study the LQR problem with a stability constraint, both in the deterministic and uncertain setting. We do not report numerical examples for the unconstrained case.\\
Before presenting the experimental results, we first recap Leslie's model to fix the notation.
\subsection{Introduction to Leslie population growth model}
The aim is to describe the dynamics of a population that can be subdivided into age classes. Each element of the state vector describes the numerosity of one of these classes. Following the notation of \eqref{eq:dyn}, matrix $G$ represents the immigration and is usually taken as the identity matrix: this implies that there are as many inputs as states. The system matrix $F$ is given by 
\begin{equation*}
F = \begin{bmatrix}
\nu_1 & \nu_2 &\cdots & \nu_{n-1} &\nu_n\\
\kappa_1\\
& \kappa_2\\
&  &\ddots \\
& & &\kappa_{n-1} & 0
\end{bmatrix}
\end{equation*}
where $\nu_i \geq 0$ is the fecundity, i.e. the average number of newborns that a member of the $i$-th class expects between $t$ and $t+1$, and $\kappa_i \in (0,1)$ represents the survival rate of the $i$-th class, $i=1...n$.

\subsection{Tests on Leslie model: deterministic setup}\label{subsec:LeslieDet}
Ths first experiment compares approaches (S0), (S1), (S2) and (S$_{\infty}$) in the deterministic case. We perform a Monte Carlo test on 100 mostly unstable Leslie models with $n=m=5$ and random parameters $\nu_i$ uniform in $(0,4)$ and $\kappa_i$ uniform in $(0,1)$. The objective matrices $Q$, $R$ and $S$ in \eqref{pb} are $Q = diag(5,4,3,2,1)$, $R=5I_5$ and $S=Q$, while the  time horizon is set to $T=8$. The initial condition is taken as $\bar{x} = [5 \;0_{1,4}]^{\top}$. Parameters $\xi$ and $\mu$ entering approach (S0) are chosen $10^{-5}$ and $0.01$ respectively \footnote{These values are the same for all numerical experiments of this Section.}. For each system we compute the feedback matrices $K_0$, $K_1$, $K_2$ and $K_{\infty}$ corresponding to the four approaches. 
\paragraph*{Remark}
In all Figures concerning this experiment, we do not present the results obtained by (S2); nevertheless, the performance of this method is described in the text below.\\

\begin{figure}[!htbp]
\centering
\includegraphics[scale=0.5]{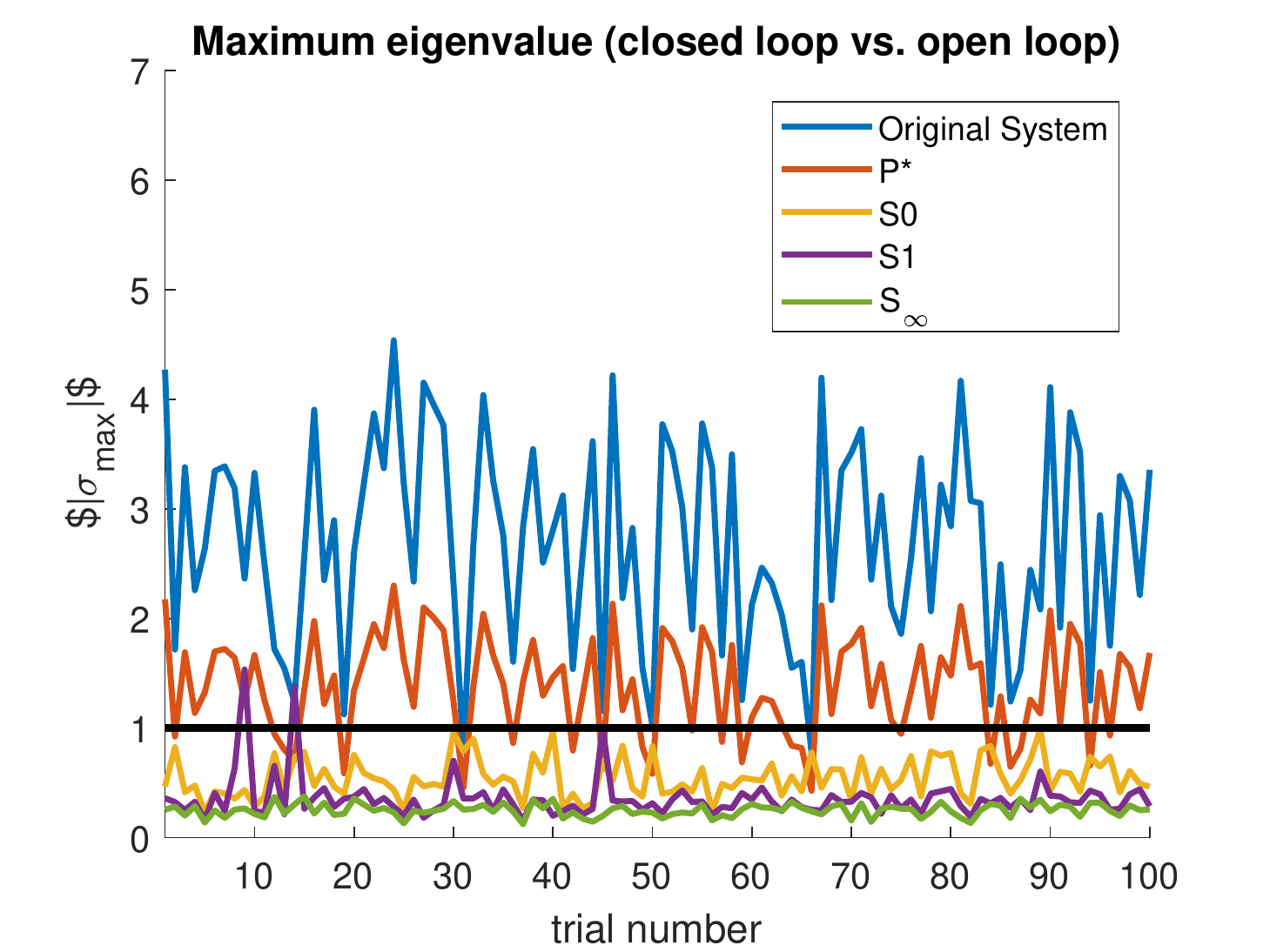}
\caption{Maximum eigenvalue}
\label{fig:maxEig}
\end{figure}

In Figure \ref{fig:maxEig} we study the performance in terms of stability: to do this, we consider the eigenvalue yielding the maximum absolute value, $|\sigma_{max}|$. In particular, we compare the ones yielded in closed-loop by $K_{0,1,2,\infty}$ with the open-loop one, and the one yielded by the optimal LQR solution $P^*$. The straight black line in the figure denotes the threshold between stable and unstable eigenvalues. From Figure \ref{fig:maxEig} it is possible to see that the original system, as well as the closed-loop obtained by $P^*$, are mostly unstable. It can be noticed that the most stabilizing one is (S$_{\infty}$). It is followed by (S1), which in fact tries to find an approximation of the stabilizing solution of the ARE; nevertheless, we have to notice that some systems were not successfully stabilized. Then, it can be noticed that (S0) is the one that shows the most robust performance.

\begin{figure}[!htbp]
\centering
\includegraphics[scale=0.5]{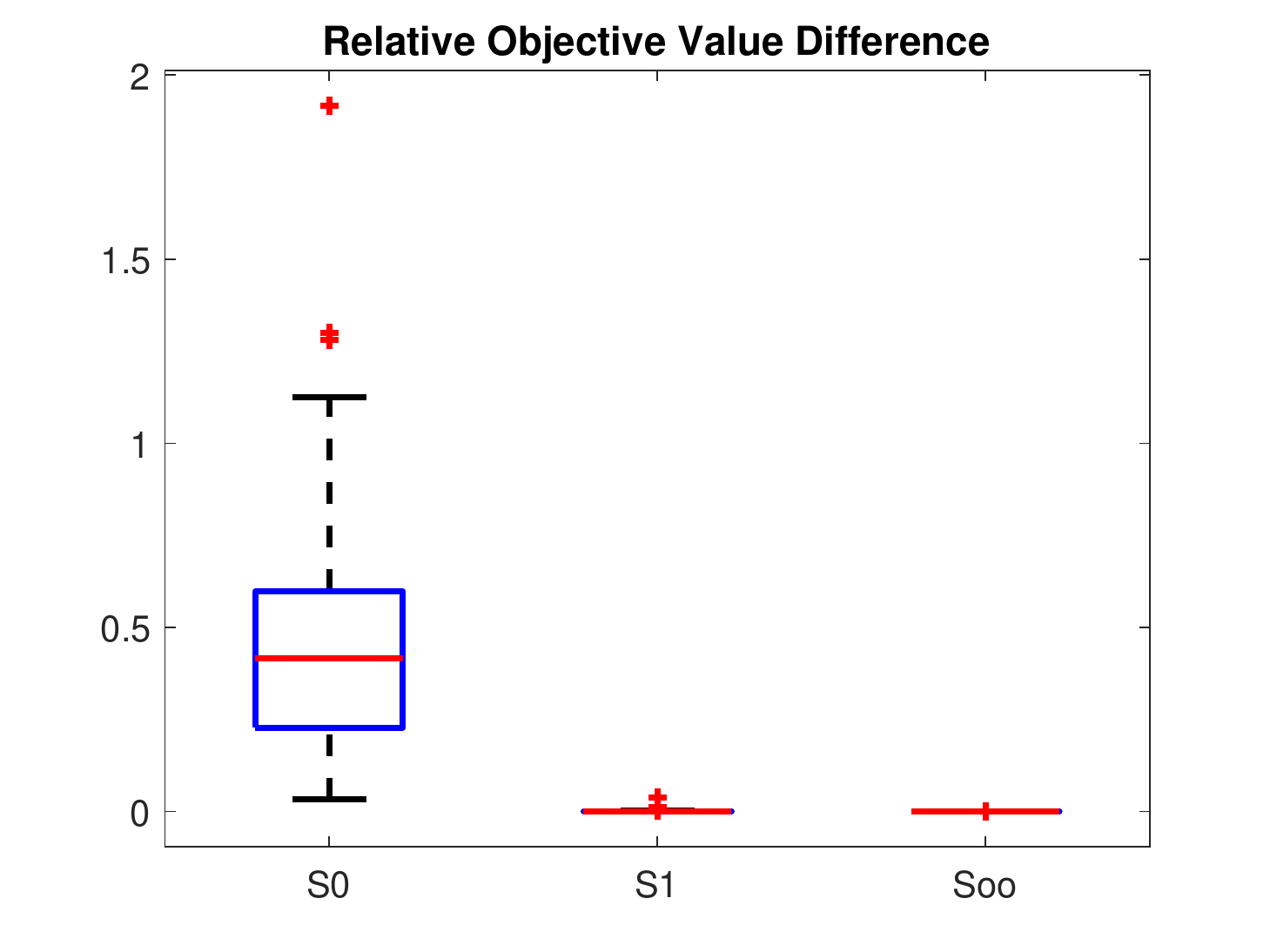}
\caption{Objective value comparison: $\frac{J_{\square}-J^*}{J^*}$.} 
\label{fig:objectiveValueDifference}
\end{figure}

We now  compare the performance of these approaches with respect to minimizing the optimal control cost. We know that the minimum is achieved by $P^{*}$, which is the solution of the classic LQR problem discussed in Section \ref{sec:classic}. The performance measure we consider is the relative objective value difference, i.e. $(J_{\square} - J^{*})/J^*$, where $J^*$ is the minimum and $J_{\square}$ is the cost yielded by approaches $\square = 0,1,2,\infty$. The resulting boxplots are shown in Figure \ref{fig:objectiveValueDifference}. Approach (S2) is largely outperformed by the others and its boxplot is therefore omitted.  (S0) yields acceptable results, but it is outperformed by (S1); on the other hand, the first guaranteed robust stability on each random system. Finally, (S$_{\infty}$) performs similarly to (S1) in terms of cost value. \\ 

\begin{figure}[!htbp]
\centering
\includegraphics[scale=0.5]{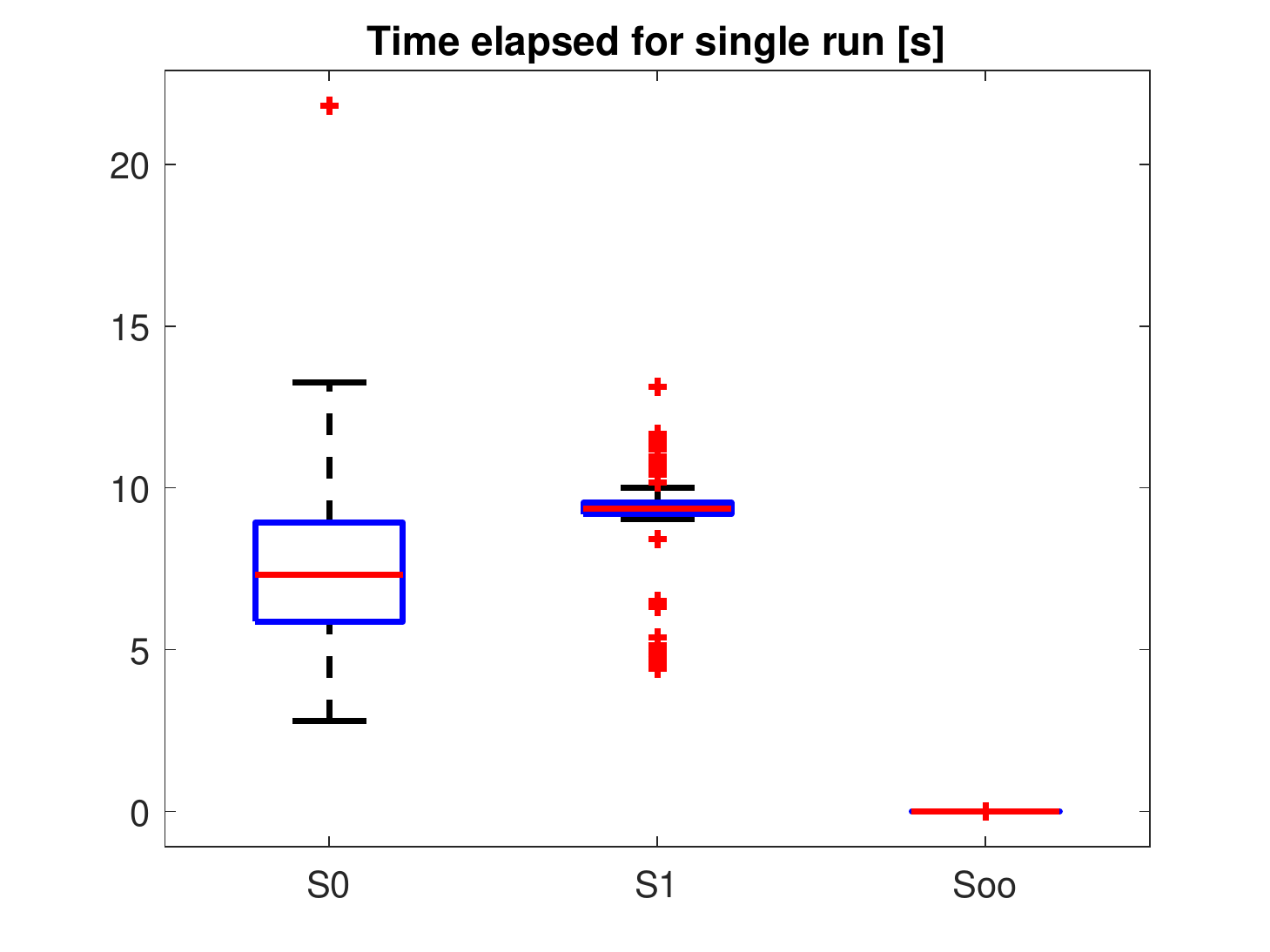}
\caption{Computation time. }
\label{fig:elapsedTime}
\end{figure}

To compare computational load, Figure \ref{fig:elapsedTime} presents the boxplot of the 100 runs.  (S$_{\infty}$) is the fastest, because its computation does not depend on $T$. Approach (S0) shows large variability in the running time, but overall has a better performance w.r.t. (S1). Moreover, (S2) is (S1) because the latter is more nonlinear, 
and therefore Nelder-Mead's algorithm is slower in contracting/reflecting to find a function minimum.

\subsection{Test with a simple non-detectable system}
In the previous tests, (S$_\infty$) seemed to be the best approach: it stabilizes effectively the system and yields in a fast way an objective value comparable with the optimal one. Nevertheless, when the detectability hypothesis in Theorem \ref{thm:stabARE} does not hold, then (S$_{\infty}$) cannot be applied, since numerical algorithms provide the output only if a stabilizing solution for the ARE exists. If this is the case, we can see from the following experiment that all other approaches are able to provide a stabilizing solution.

Let us consider $F = diag(2,1)$, $G = 1_2$, $R = 1$, $Q = diag(1,0)$, $S=I_2$ and $\bar{x} = [1\,0]^{\top}$. It can be readily checked that $(F,Q^{1/2})$ is not detectable. We apply methods (S0), (S1), (S2): in this case, (S$_{\infty}$) concides with the open loop. 

\begin{figure}[!htbp]
\centering
\includegraphics[scale=0.5]{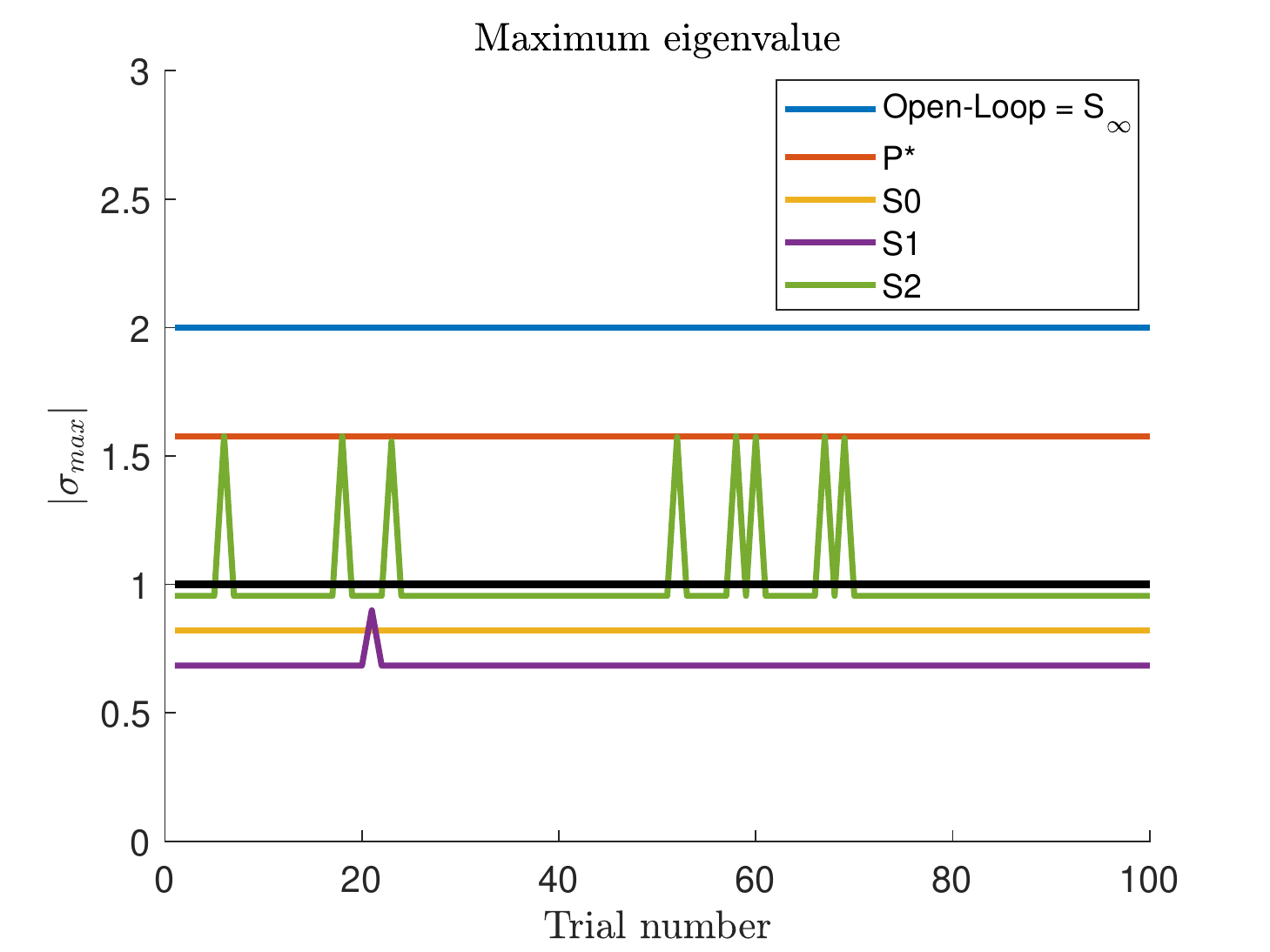}
\caption{Maximum eigenvalue in a case in which $(F,Q^{1/2})$ is not detectable}
\label{fig:nondetec}
\end{figure}

The experiment is run 100 times on the same deterministic system, to highlight the dependence of all approaches from the initial conditions, that are set random. We can notice in particular that (S2) for 8/100 times minimizes its objective by putting $LL^{\top}=S$, thus yielding solution $P^{*}$ which is optimal in terms of the cost but is not stabilizing. 
Approaches (S0) and (S1)  are reliable in terms of stability, though (S1) is subject to oscillations. \\
In Table \ref{tab:1} we compare the approaches in terms of the relative cost difference.

\begin{table}[ht]
\centering
\begin{tabular}{|c|c|}
\hline 
Approach & $(J_{\square} - J^{*})/J^*$\\
\hline
\hline
(S0) & 0.1843\\
(S1) & 0.0366\\
(S2) & 37.1981\\
(S$_\infty$)=open loop &1.5$\times10^4$\\
\hline
\end{tabular}
\caption{Comparison of the objective values with metric $\frac{J_{\square}-J^*}{J^*}$}
\label{tab:1}
\end{table}
Again, (S0) performs slightly worse than (S1), that in this case is always able to stabilize the system.

\subsection{Test on Leslie model: scenario setup}
We now test approaches (S0), (S1), (S2) and (S$_{\infty}$) with the scenario approach on the following Leslie model of dimension $n=5$ presenting parametric uncertainty:
\begin{equation}\label{experimentScenario}
\begin{cases}
\nu_1=1.11 \\
\nu_2=2.05\\
\nu_3=1.79\\
\nu_4=2.37\\
\nu_5=1.10
\end{cases}\qquad \begin{cases}
\kappa_1=0.97+\delta_1\\
\kappa_2=0.86+\delta_2\\
\kappa_3=0.37+\delta_3\\
\kappa_4=0.09+\delta_4
 \end{cases}
\end{equation}
where $\delta_i$ is distributed as a uniform random variable over $[-0.4, 0.4]$. This system turns out to be always unstable: we aim at finding a stabilizing feedback fixing parameter $\beta=0.05$ and having $N=50$ scenarios available from past measurements. \\
The objective matrices $Q$, $R$ and $S$ in \eqref{pb} are chosen as in Section \ref{subsec:LeslieDet}, i.e. $Q = diag(5,4,3,2,1)$, $R=5I_5$ and $S=Q$, while the considered time horizon is set to $T=8$. As initial value, we consider $\bar{x}$ as in the previous experiment. 
We run all the four approaches on these data and test the resulting feedback matrices on 100 new samples of the system, and \textit{all approaches yield stabilizing solutions.} \\
For robustness guarantees, the min-max structure of the problem yields a support constraint set of cardinality one: hence, according to \eqref{eq:epsilon},
\begin{displaymath}
\epsilon(s_{N,\square}^*) = \epsilon(1) = 0.1981 \qquad \text{ for } \square = 0, 1, 2, \infty.
\end{displaymath}
Objective value comparisons and running times  are consistent with the experiment of Section \ref{subsec:LeslieDet} and are therefore omitted.

\section{Conclusions}\label{sec:conclusions}
In this paper, we have proposed a new perspective on the concept of stability in the finite-horizon LQR scenario. Our grand strategy consists of formulating an optimization problem that balances performance (in terms of cost value) and closed-loop stability. The stability-constrained problem has been also treated in a robust framework via the scenario approach.\\
First, we stated the classical LQR problem in four equivalent ways denoted by (P1), (P2), (P3) and (P4). Then, generalizing (P1), the solution of the unconstrained robust LQR problem has been found. The resulting problem is convex and, hence, it provides a priori guarantees on robustness.
The other three formulations (P2), (P3) and (P4) are exploited to deal with the stability constraint, first in the deterministic and then in the robust case. 
Their generalizations lead to the following approaches:
\begin{itemize}
\item \textbf{(S0)} uses the formulation (P2) and includes additional variables in order to include the stability constraint as a LMI inspired by Lyapunov's Theorem. The overall problem is solved by alternating 
minimization over two structured convex subproblems. It yields a  good trade-off between cost, stability and computational time. 
\item \textbf{(S$_{\infty}$)} relies on the theory concerning the Algebraic Riccati Equation, which yields the optimal solution of the infinite-horizon LQR problem. This approach performs best w.r.t. cost, stability and computational time, but fails whenever the detectability condition in Theorem \ref{thm:stabARE} is not satisfied. 
\item \textbf{(S1) and (S2)} derive from (P4) and (P3) respectively. Both require solving nonconvex unconstrained optimization problems that minimize the 2-norm of the state. Moreover, in view of the structure assigned to the input, their output can be regarded as the approximator of the stabilizing solution of the ARE. From this point of view, (S1) is very effective but its drawback is the computational time due to its nonlinearity. On the other hand, (S2) is solved faster but is less reliable in terms of stability and cost. Both can be used when (S$_{\infty}$) cannot be applied due to non-detectability.
\end{itemize}

All of these methods involve nonconvex programs. Nevertheless, convergence analyses and numerical results show that they can be efficiently implemented, obtaining a good trade-off between cost, stability and computational time. When uncertainty is treated using the scenario approach, the robustness of the solution can be effectively quantified a posteriori with a suitable support constraint selection.\\ 
In the past control, system identification and statistical literature, the use of nonconvex procedures and models was  looked upon with suspicion. More recently, the utility of such techniques has become evident e.g. in the machine learning field with the adoption of deep networks \cite{GoodfellowBengio}. Along this line, this paper shows that nonconvex approaches can be effective and computationally efficient also in the optimal control field. To make them even more effective, further research developments may concern (S1) and (S2). For instance, by exploring gradient-based methods and Algorithmic Differentiation \cite{bradjim} it could be possible to control also more complex large-scale and interconnected systems.

\bibliographystyle{plain}        
\bibliography{mybibfile}           

\begin{thebibliography}{10}

\bibitem{alamotempocamacho}
T.~Alamo, R.~Tempo, and E.~Camacho.
\newblock Randomized strategies for probabilistic solutions of uncertain
  feasibility and optimization problems.
\newblock {\em Automatic Control, IEEE Transactions on}, 54:2545 -- 2559, 12
  2009.

\bibitem{Anderson}
B.~D.~O. Anderson and J.~B. Moore.
\newblock {\em Optimal Control: Linear Quadratic Methods}.
\newblock Prentice-Hall, Inc., Upper Saddle River, NJ, USA, 1990.

\bibitem{aravkin2012estimating}
Aleksandr~Y Aravkin and Tristan Van~Leeuwen.
\newblock Estimating nuisance parameters in inverse problems.
\newblock {\em Inverse Problems}, 28(11):115016, 2012.

\bibitem{ArnoldLaub}
W.~F. {Arnold} and A.~J. {Laub}.
\newblock Generalized eigenproblem algorithms and software for algebraic
  riccati equations.
\newblock {\em Proceedings of the IEEE}, 72(12):1746--1754, Dec 1984.

\bibitem{bradjim}
B.M. Bell and J.~Burke.
\newblock {\em Algorithmic Differentiation of Implicit Functions and Optimal
  Values}, volume~64, pages 67--77.
\newblock 08 2008.

\bibitem{Bellman}
R.~E. Bellman.
\newblock {\em Dynamic Programming}.
\newblock Dover Publications, Inc., New York, NY, USA, 2003.

\bibitem{BEMPORAD19981255}
A.~Bemporad.
\newblock A predictive controller with artificial lyapunov function for linear
  systems with input/state constraints.

\bibitem{Bertsekas}
D.~P. Bertsekas.
\newblock {\em Dynamic Programming and Optimal Control}.
\newblock Athena Scientific, 2nd edition, 2000.

\bibitem{Bitmead1991RiccatiDA}
R.~R. Bitmead and M.~Gevers.
\newblock Riccati difference and differential equations: Convergence,
  monotonicity and stability.
\newblock 1991.

\bibitem{BEFB:94}
S.~Boyd, L.~{El~{G}haoui}, E.~Feron, and V.~Balakrishnan.
\newblock {\em Linear Matrix Inequalities in System and Control Theory},
  volume~15 of {\em Studies in Applied Mathematics}.
\newblock {SIAM}, Philadelphia, PA, June 1994.

\bibitem{Boyd}
S.~Boyd and L.~Vandenberghe.
\newblock {\em Convex Optimization}.
\newblock Cambridge University Press, New York, NY, USA, 2004.

\bibitem{BressanPiccoli}
A.~Bressan and B.~Piccoli.
\newblock Introduction to the mathematical theory of control.
\newblock 01 2007.

\bibitem{brysonho}
A.~E. Bryson and Y.~Ho.
\newblock {\em Applied optimal control: optimization, estimation, and control}.
\newblock Hemisphere Pub. Corp. ; distributed by Halsted Press Washington : New
  York, rev. printing. edition, 1975.

\bibitem{Calafiore2013RobustMP}
G.~Calafiore and L.~Fagiano.
\newblock Robust model predictive control via scenario optimization.
\newblock {\em IEEE Transactions on Automatic Control}, 58:219--224, 2013.

\bibitem{calafiorecampi}
G.~C. {Calafiore} and M.~C. {Campi}.
\newblock The scenario approach to robust control design.
\newblock {\em IEEE Transactions on Automatic Control}, 51(5):742--753, May
  2006.

\bibitem{CaE:14}
G.C. Calafiore and L.~{El Ghaoui}.
\newblock {\em Optimization Models}.
\newblock Control systems and optimization series. Cambridge University Press,
  October 2014.

\bibitem{CampiGarattiM03}
M.C. Campi and S.~Garatti.
\newblock {\em Introduction to the Scenario Approach}.
\newblock Society for Industrial and Applied Mathematics, Philadelphia, PA,
  2018.

\bibitem{campinonconv}
M.C. Campi, S.~Garatti, and F.~Ramponi.
\newblock A general scenario theory for nonconvex optimization and decision
  making.
\newblock {\em IEEE Transactions on Automatic Control}, PP:1--1, 02 2018.

\bibitem{denicolaomagnistabilitymmpc}
G.~De~Nicolao, L.~Magni, and R.~Scattolini.
\newblock Stability and robustness of nonlinear receding horizon control.
\newblock In Frank Allg{\"o}wer and Alex Zheng, editors, {\em Nonlinear Model
  Predictive Control}, pages "3--22", Basel, 2000. Birkh{\"a}user Basel.

\bibitem{DEOLIVEIRA1999}
M.C. de~Oliveira, J.~Bernussou, and J.C. Geromel.
\newblock A new discrete-time robust stability condition.
\newblock {\em Systems and Control Letters}, 37(4):261 -- 265, 1999.

\bibitem{finsler1937uber}
P.~Finsler.
\newblock {\em Uber das Vorkommen definiter und semidefiniter Formen in Scharen
  quadratischer Formen}.
\newblock In aedibus O. Fussli, 1937.

\bibitem{GoodfellowBengio}
I.~Goodfellow, Y.~Bengio, and A.~Courville.
\newblock {\em Deep Learning}.
\newblock The MIT Press, 2016.

\bibitem{grammaticononconv}
S.~Grammatico, X.~Zhang, K.~Margellos, P.~Goulart, and J.~Lygeros.
\newblock A scenario approach for non-convex control design.
\newblock {\em IEEE Transactions on Automatic Control}, 61, 02 2016.

\bibitem{grune2017nonlinear}
L.~Gr{\"u}ne and J.~Pannek.
\newblock Nonlinear model predictive control.
\newblock In {\em Nonlinear Model Predictive Control}, pages 45--69. Springer,
  2017.

\bibitem{Helly1923}
E.~Helly.
\newblock \"uber mengen konvexer k\"orper mit gemeinschaftlichen punkte.
\newblock {\em Jahresbericht der Deutschen Mathematiker-Vereinigung},
  32:175--176, 1923.

\bibitem{hou2015dynamic}
Chengjun Hou.
\newblock {\em Dynamic programming under parametric uncertainty with
  applications in cyber security and project management}.
\newblock PhD thesis, The Ohio State University, 2015.

\bibitem{nonconvOptMachineLearning}
P.~{Jain} and P.~{Kar}.
\newblock {\em Non-convex Optimization for Machine Learning}.
\newblock 2017.

\bibitem{Kalman60contributionsto}
R.E. Kalman.
\newblock Contributions to the theory of optimal control, 1960.

\bibitem{kleinman}
D.~{Kleinman}.
\newblock An easy way to stabilize a linear constant system.
\newblock {\em IEEE Transactions on Automatic Control}, 15(6):692--692,
  December 1970.

\bibitem{kuhn1951}
H.~W. Kuhn and A.~W. Tucker.
\newblock Nonlinear programming.
\newblock In {\em Proceedings of the Second Berkeley Symposium on Mathematical
  Statistics and Probability}, pages 481--492, Berkeley, Calif., 1951.
  University of California Press.

\bibitem{kwakernaak}
H.~Kwakernaak and R.~Sivan.
\newblock {\em Linear optimal control systems}.
\newblock Wiley Interscience New York, 1972.

\bibitem{kwon}
W.H. Kwon, A.M. Bruckstein, and T.~Kailath.
\newblock Stabilizing state-feedback design via the moving horizon method.
\newblock {\em International Journal of Control}, 37(3):631--643, 1983.

\bibitem{1961iii}
G.~La~Salle and S.~Lefschetz.
\newblock {\em Stability by Liapunov's Direct Method with Applications},
  volume~4 of {\em Mathematics in Science and Engineering}.
\newblock Elsevier, 1961.

\bibitem{lancaster1995algebraic}
P.~Lancaster and L.~Rodman.
\newblock {\em Algebraic Riccati Equations}.
\newblock Oxford science publications. Clarendon Press, 1995.

\bibitem{leslie}
P.~H. Leslie.
\newblock On the use of matrices in certain population mathematics.
\newblock {\em Biometrika}, 33(3):183--212, 1945.

\bibitem{lewis}
F.L. Lewis.
\newblock {\em Optimal control}.
\newblock A Wiley-Interscience publication. Wiley, 1986.

\bibitem{LuenbergerVectorSpaceOpt}
D.~G. Luenberger.
\newblock {\em Optimization by Vector Space Methods}.
\newblock John Wiley \& Sons, Inc., New York, NY, USA, 1st edition, 1997.

\bibitem{luenberger1979introduction}
D.G. Luenberger.
\newblock {\em Introduction to Dynamic Systems: Theory, Models, and
  Applications}.
\newblock Wiley, 1979.

\bibitem{Magnus1985MatrixDC}
J.~R. Magnus and H.~Neudecker.
\newblock Matrix differential calculus with applications to simple, hadamard,
  and kronecker products.
\newblock 1985.

\bibitem{martenssonrantzer}
K.~M{\aa}rtensson and Anders Rantzer.
\newblock Synthesis of structured controllers for large-scale systems.
\newblock 07 2019.

\bibitem{Mayne2000ConstrainedMP}
D.Q Mayne, J.B Rawlings, C.V. Rao, and P.O.M Scokaert.
\newblock Constrained model predictive control: Stability and optimality.
\newblock {\em Automatica}, 36:789--814, 2000.

\bibitem{NoceWrig06}
J.~Nocedal and S.~J. Wright.
\newblock {\em Numerical Optimization}.
\newblock Springer, New York, NY, USA, second edition, 2006.

\bibitem{pontryagin}
L.S. Pontryagin.
\newblock {\em Mathematical Theory of Optimal Processes}.
\newblock Classics of Soviet Mathematics. Taylor \& Francis, 1987.

\bibitem{mpc}
S.~V. Raković and W.~S.~Levine.
\newblock {\em Handbook of Model Predictive Control}.
\newblock 09 2018.

\bibitem{rawlings2009model}
J.B. Rawlings and D.Q. Mayne.
\newblock {\em Model Predictive Control: Theory and Design}.
\newblock Nob Hill Pub., 2009.

\bibitem{recht}
B.~Recht.
\newblock A tour of reinforcement learning: The view from continuous control.
\newblock {\em Annual Review of Control, Robotics, and Autonomous Systems},
  2(1):253--279, 2019.

\bibitem{thomas}
Y.~A. {Thomas}.
\newblock Linear quadratic optimal estimation and control with receding
  horizon.
\newblock {\em Electronics Letters}, 11(1):19--21, January 1975.

\bibitem{todorovconf}
E.~Todorov.
\newblock General duality between optimal control and estimation.
\newblock volume~47, pages 4286 -- 4292, 01 2009.

\bibitem{TsengAltMin}
P.~Tseng.
\newblock Applications of a splitting algorithm to decomposition in convex
  programming and variational inequalities.
\newblock {\em SIAM Journal on Control and Optimization}, 29(1):119--138, 1991.

\bibitem{tseng2001convergence}
P.~Tseng.
\newblock Convergence of a block coordinate descent method for
  nondifferentiable minimization.
\newblock {\em Journal of optimization theory and applications},
  109(3):475--494, 2001.

\bibitem{Vapnik1998}
V.N. Vapnik.
\newblock {\em Statistical Learning Theory}.
\newblock Wiley-Interscience, 1998.

\end{thebibliography}

\appendix

\section{On the derivation of expression \eqref{eq:gradient}}
\subsection{Catalogue of useful properties}
Before performing the computation for our case, we first recall some useful properties of differentials, vec/Kronecker/trace operators and commutation matrices \cite{Magnus1985MatrixDC} . We will denote with $X$, $Y$, $W$ and $Z$ general matrices of suitable dimensions. 
\begin{enumerate}[label=\subscript{r}{{\arabic*}})]
\item \label{r1} $vec(YXZ)=(Z^{\top}\otimes Y)vec(X)$;
\item \label{r2} $(X \otimes Y)^{\top} = X^{\top}\otimes Y^{\top}$;
\item \label{r3} $(X\otimes Y)(W\otimes Z) = (XW)
\otimes(YZ)$;
\item \label{r4} $vec(dX) = dvec(X)$;
\item \label{r5} $(dX^{\top}) = (dX)^{\top}$;
\item \label{r6} $tr(X+Y) = tr(X) + tr(Y)$;
\item \label{r7} $tr(XY)=tr(YX)$;
\item \label{r8} $tr(X^{\top}) = tr(X)$;
\item \label{r9} $tr(Z^{\top}dX) = (vec(Z))^{\top}vec(dX)$;
\item \label{r10} if $X \in \mathbb{R}^{d_1\times d_2}$, then $vec(X^{\top})=\mathscr{K}_{d_1d_2}vec(X)$. Notice that $\mathscr{K}_{d_1d_2}$ does not depend on $X$, but only on its dimension;
\item \label{r11} $\mathscr{K}_{d_1d_2}^{-1} = \mathscr{K}_{d_1d_2}^{\top}$;
\item \label{r12} considering $X \in \mathbb{R}^{d_1\times d_2}$ and $Y \in \mathbb{R}^{d_3 \times d_4}$, then their Kronecker product is a $d_1d_3 \times d_2d_4$ matrix such that $\mathscr{K}_{d_1d_3}(X \otimes Y) = (Y \otimes X)\mathscr{K}_{d_2d_4}$;
\item \label{r13} $\mathscr{K}_{d_1^2} = \mathscr{K}_{d_1^2}^{\top}$.
\end{enumerate}

The computation of the gradient of a certain function $\phi:\mathbb{R}^{d_1\times d_2} \rightarrow \mathbb{R}^{d_3\times d_4}$ w.r.t. $X \in \mathbb{R}^{d_1\times d_2}$ consists in the following steps:
\begin{enumerate}[label=\subscript{v}{{\arabic*}})]
\item compute the differential $d\phi$ and use properties \eqref{r1}...\eqref{r13} to express it in terms of $vec(dX)$, possibly isolating it as the rightmost factor;
\item exploiting \eqref{r4}, compute $d\phi / dvec(X)$;
\item use properties \eqref{r1}...\eqref{r13} to obtain an expression that depends on $X$ and not $vec(X)$.
\end{enumerate}

\subsection{Computations}
We recall that the objective is $$\bar{\mathcal{L}}(x,\lambda,K) + \frac{1}{2\mu}\|KC-D\|_F^2=\Omega(K),$$
where $\bar{\mathcal{L}}(x,\lambda,K)$ is defined as in \eqref{eq:S0obj}. The variable of interest will be $K$ only, as discussed in Section \ref{sec:S0}.
\paragraph*{Step (v1)} We first aim at finding $d\Omega$: we consider its two terms separately.
As regards $d\bar{\mathcal{L}}(x,\lambda,K)$, we first write the quadratic terms by means of 2-norms and get 
\begin{align*}
&d\bar{\mathcal{L}}(x,\lambda,K) = \sum_{t=0}^{T-1} \| R^{1/2}(dK)x_t \|^2 - \lambda_{t+1}^{\top}G(dK)x_t\\
&= \sum_{t=0}^{T-1} \| (x_t^{\top} \otimes R^{1/2})vec(dK)\|^2 -\\ &\qquad (x_t^{\top} \otimes \lambda_{t+1}^{\top}G)vec(dK) \qquad \text{(by \eqref{r1}}\\
&= \sum_{t=0}^{T-1} (vec(dK))^{\top}(x_t \otimes R^{1/2})(x_t^{\top}\otimes R^{1/2})vec(dK) -\\ &\qquad - (x_t^{\top} \otimes \lambda_{t+1}^{\top}G)vec(dK) \qquad \text{(by \eqref{r2}}\\
&= \sum_{t=0}^{T-1}  (vec(dK))^{\top}(x_tx_t^{\top} \otimes R)vec(dK) - \\ & \qquad -(x_t^{\top} \otimes \lambda_{t+1}^{\top}G)vec(dK) \qquad \text{(by \eqref{r3}}.\\
\end{align*}
As regards $\|KC-D\|_F^2$, we first recall that for any matrix $Y$ it holds that $\|Y\|_F^2 = tr(Y^{\top}Y)$. Then, the differential reads 
\begin{align*}
&dtr((KC-D)^{\top}(KC-D)) = \\& \underbrace{tr(C^{\top}(dK^{\top})KC) + tr(C^{\top}K^{\top}(dK)C)}_{(a)} -\\ &-\underbrace{tr(D^{\top}(dK)C)}_{(b)} - \underbrace{tr(C^{\top}(dK^{\top})D)}_{(c)}
\end{align*}
From simple computations it follows that
\begin{align*}
(a)&=2tr(CC^{\top}K^{\top}(dK)) \qquad \text{ (by (\ref{r7}, \eqref{r8}}\\
&=2(vec(KCC^{\top}))^{\top}vec(dK) \qquad \text{ (by \eqref{r9}}
\end{align*}
\begin{align*}
(b)&=tr(CD^{\top}(dK)) \qquad\text{ (by \eqref{r7}}\\
&=vec(CD^{\top}))^{\top}vec(dK) \qquad\text{ (by \eqref{r9}}
\end{align*}
\begin{align*}
(c)&=tr(DC^{\top}(dK^{\top})) \quad \text{(by \eqref{r7}}\\
&= (vec(CD^{\top}))^{\top}vec(dK^{\top}) \quad \text{(by \eqref{r9}}\\
&= ((D \otimes C)vec(I_n))^{\top}\mathscr{K}_{mn}vec(dK) \quad \text{(by (\ref{r1} and \eqref{r10}}\\
&= (vec(I_n))^{\top}(D^{\top} \otimes C^{\top})\mathscr{K}_{mn}vec(dK) \quad \text{(by \eqref{r2}}\\
&= (vec(I_n))^{\top}\mathscr{K}_{n^2}(C^{\top} \otimes D^{\top})vec(dK) \quad \text{(by \eqref{r12}}\\
&= [\mathscr{K}_{n^2}vec(I_n)]^{\top}(C^{\top} \otimes D^{\top})vec(dK) \quad \text{(by (\ref{r12} and \eqref{r13}}\\
&= [(C \otimes D)vec(I_n)]^{\top}vec(dK)\\
&= (vec(DC^{\top}))^{\top}vec(dK \quad \text{(by \eqref{r1}}\\
&=(b).
\end{align*}
\paragraph*{Step (v2)} Using the expression for $d\Omega$ above computed, we get
\begin{align*}
\frac{d\Omega}{dvec(K)} &= \sum_{t=0}^{T-1}2(x_tx_t^{\top} \otimes R)vec(dK) - (x_t \otimes G^{\top}\lambda_{t+1}) \\&\quad + \frac{1}{2\mu}(2(vec(KCC^{\top}))^{\top} - 2(vec(DC^{\top}))^{\top}).
\end{align*}

\paragraph*{Step (v3)} We finally get \eqref{eq:gradient} by exploiting (\ref{r1} and (\ref{r9}.
The values from $x_t$ and $\lambda_t$ can be recovered from their matrix form of \eqref{eq:C1optcond}. If we call $F + GK =: F_K$, then the expressions for $x_t$ and $\lambda_t$ are $x_t = F_K^t\bar{x}$ and 
\begin{displaymath}
\scriptstyle{\lambda_t = -2\Big[(F_K^{\top})^{T-t}SF_K^T+\sum_{i=0}^{T-t-1}(F_K^{\top})^i(K^{\top}RK + Q)F_K^{t+i}\Big]\bar{x}} 
\end{displaymath}
respectively.


\end{document}